`%hodes.tex: 
%%a Plain TeX file by Shalosh B. Ekhad and Doron Zeilberger (5 pages)

%begin macros

\baselineskip=14pt
\parskip=10pt
\def\halmos{\hbox{\vrule height0.15cm width0.01cm\vbox{\hrule height
  0.01cm width0.2cm \vskip0.15cm \hrule height 0.01cm width0.2cm}\vrule
  height0.15cm width 0.01cm}}
\font\eightrm=cmr8 
\font\eighttt=cmtt8
\magnification=\magstephalf

\def\1{{\overline{1}}}
\def\2{{\overline{2}}}
\parindent=0pt
\overfullrule=0in

\def\frac#1#2{{#1 \over #2}}
%\headline={\rm  \ifodd\pageno  \RightHead  \else  \LeftHead  \fi}
%\def\RightHead{\centerline{
%Title
%}}
%\def\LeftHead{ \centerline{Doron Zeilberger}}
%end macros
\centerline
{\bf Every Fifth\footnote{$^1$}{\eightrm  More precisely: \hfill \break
every  $5.000000000000000000000000000000000000000000000000000000000000005415558$-th real number}
Real Number is Evil }
\smallskip
\centerline
{\it Doron ZEILBERGER}

\smallskip

{\it In fond memory of my beloved cousin Matti Weiss (1945-2025) who ignited my love of mathematics}

\bigskip
{\bf Abstract:} Fifteen years ago, then-Carleton-undergrad, Isaac Hodes, proved that the Golden Ratio is evil. In this modest contribution to human knowledge,
we show that in fact, every fifth real number is evil, and we present lots of other interesting numbers that are evil. We also show
(in addition to many other fascinating facts), that the expected evil-location of a random evil number is $148.185185185185185\dots$ .
It follows that the Golden Ratio is a fairly average evil real number, since, as first computed by Hodes, its evil location is $146$.

\bigskip

{\bf Maple package}: This article is accompanied by a Maple package {\tt Hodes.txt}, available from

{\tt http://sites.math.rutgers.edu/\~{}zeilberg/tokhniot/Hodes.txt } \quad .

The front of this article

{\tt https://sites.math.rutgers.edu/\~{}zeilberg/mamarim/mamarimhtml/hodes.html } \quad,

contains numerous output files, some of which will be referred to later.

{\bf Original Definition:} ([H]) A real number is {\bf evil} if the sequence of digits in its decimal representation to the right of the decimal point (in other words, its {\it fractal} part)
has the property that one of its {\it partial sums} is the satanic number $666$.

Hodes[H] showed that, according to his definition, the golden ratio $\frac{1+\sqrt{5}}{2}$, is evil, and the evil location is $146$, i.e. the sum of the first $146$ digits after
the decimal point equals $666$.

But Hodes' definition is unfair to the digits to the left of the decimal point. Also it is decimal-centric, dependent on the contingent fact that we humans have ten fingers.
Also there are lots of other definitions of being evil . For example mine is $359$ (the {\it gematria} value of {\it satan} ([shin, tet, nun]). Of course $13$ is a famous one.
$666$ is the Gematria value of ``Emperor Neron'' ([nun, resh, vav,nun], [kuf,samech, resh]).

A more satisfactory and general definition is the following.

{\bf My Definition of Evil Real Number}: A real number is $n$-evil to base $b$ if removing the decimal point (more generally: the base $b$ `point'), one of the partial sums of
its sequence of digits equals $n$.

Note that according to my definition, the Golden Ratio $\phi$, is not evil, but $\phi-1=\frac{\sqrt{5} -1}{2}$ is.

It turns out that the Golden Ratio $\phi$ (or according to our definition $\phi -1$) is not that special, and we have the following.

{\bf Fact}: The probability that a random real number (in base $10$) would have one of its partial sums equal to $666$ is {\bf exactly} 
$$
0.19999999999999999999999999999999999999999999999999999999999999978337773162864760552794625\dots \quad .
$$
It follows that one in every
$$
5.0000000000000000000000000000000000000000000000000000000000000054155\dots
$$
real numbers is evil. Let me explain.

{\bf Probability Generating Functions}

{\bf Theorem 1}:  Let $a_b(n)$ be the probability that a uniformly-at-random string of members of $\{0,1, \dots , b-1\}$ that definitely starts with a non-zero entry,
has one of its partial sums equal to $n$. We have the following generating function, let's call it $h_b(x)$:
$$
h_b(x) \, := \, \sum_{n=0}^{\infty} \, a_b(n)\, x^n =
1+ \frac{x(1-x^{b-1})}{(b-1)(1-x)}+
\frac{\left(-1+x^{b -1}\right)^{2} x^{2}}{\left(-b x +x^{b}+b -1\right) \left(b -1\right) \left(1-x \right)} \quad .
$$

In fact we have more generally:

{\bf Theorem 2}: Let $A_b(n,k)$ be the probability that a random real number in base $b$ has its $k$-th partial sum equal to $n$ (for the first  time). We have the
following bivariate generating function, let's call it $H_b(x,t)$:
$$
H_b(x,t) \, := \,
\sum_{n=0}^{\infty} \sum_{k=0}^{\infty} A_b(n,k) \, x^n \, t^k \,
= \, 1+ \frac{x(1-x^{b-1})}{(b-1)(1-x)} \, t + \frac{\left(1-x^{b -1}\right)^{2} x^{2} t^{2}}{\left(t \,x^{b}-b x +b -t \right) \left(b -1\right) \left(1-x \right)}  \quad .
$$

{\bf Proof}: The probability generating function of a single digit (assuming a fair die) is
$$
\frac{1}{b} \cdot \sum_{i=0}^{b-1} x^i \, = \, \frac{1-x^b}{b(1-x)} \quad .
$$

By fiat, the first digit is non-zero, so the probability generating function of the very first digit is:
$$
\frac{x+x^2+ \dots +x^{b-1}}{b-1} = \frac{x(1-x^{b-1})}{(b-1)(1-x)} \quad.
$$

If $k$ is the location where, for the {\bf first} time, the partial sum of the (uniformly-at-random) string of members from $\{0,1, \dots, b-1\}$ equals $n$,
that digit must be {\bf non-zero}. But a-priori $0$ is allowed, so the probability generating function for that $k^{th}$ digit  is:

$$
\frac{x+x^2+ \dots +x^{b-1}}{b} = \frac{x(1-x^{b-1})}{b \,(1-x)} \quad.
$$

Hence the probability generating function that the $k$-th partial sum would equal to $n$ for the first time is
$$
\frac{x(1-x^{b-1})}{(b-1)\,(1-x)} \cdot \left( \frac{1-x^b}{b(1-x)} \right )^{k-2}  \frac{x(1-x^{b-1})}{b (1-x)} \quad .
$$
Hence the desired bivariate probability generating function is
$$
H_b(x,t)=
1+ \frac{x(1-x^{b-1})}{(b-1)(1-x)} \, t \,+\,
\sum_{k=2}^{\infty} \, \frac{x(1-x^{b-1})}{(b-1)\,(1-x)} \cdot \left( \frac{1-x^b}{b(1-x)} \right )^{k-2}  \frac{x(1-x^{b-1})}{b \, (1-x)} \, t^k \quad 
$$
$$
=1+ \frac{x(1-x^{b-1})}{(b-1)(1-x)} \, t \,+\,
t^2\, \frac{x(1-x^{b-1})}{(b-1) \, (1-x)} \left ( 1-  \frac{1-x^b}{b(1-x)} \, t \right )^{-1}  \frac{x(1-x^{b-1})}{b\,(1-x)}
$$
$$
= \, 1+ \frac{x(1-x^{b-1})}{(b-1)(1-x)} \, t + \frac{\left(1-x^{b -1}\right)^{2} x^{2} t^{2}}{\left(t \,x^{b}-b x +b -t \right) \left(b -1\right) \left(1-x \right)} \quad . \halmos
$$

Plugging-in $t=1$ in Theorem $2$ gives Theorem $1$.

{\eighttt This is implemented in procedures {\tt Gxr(x,r)} and {\tt Gxtr(x,t,r)} in the Maple package {\tt Hodes.txt}} \quad .

In particular, for base $10$, we have that the probability that a random real number would have one of its partial sums equal to $666$ is the
coefficient of $x^{666}$ in the Maclaurin expansion of
$$
-\frac{9}{x^{9}+x^{8}+x^{7}+x^{6}+x^{5}+x^{4}+x^{3}+x^{2}+x -9} \quad,
$$
that happens to be the number mentioned above.

{\bf Automatic Generation of moments}

We are interested in the random variable `evil-location' in the sample space of all $n$-evil numbers in base $b$. Using partial fractions and
Maple, one can obtain very precise asymptotic expressions not only for the (conditional on being evil) expectation, and variance, but as
many moments as one desires (we went up to $16$-th, but one can easily go further).
In order to get the generating function of the $i$-th moment, Maple computed, $(t\frac{d}{dt})^i H_b(x,t) \vert_{t=1}$, a certain rational function of $x$ whose
denominator can be factored as $(1-x)^{i+1} Q(x)$ where all the roots of $Q(x)$ have absolute value larger than $1$. Then we (or rather Maple)
converted it to {\it partial fractions} using {\tt convert(f,parfrac)}, and then extracted the coefficient of $x^n$ in the
Maclaurin expansion, ignoring the exponentially decaying terms. Once we had
the moments, we got the {\it central moments}, $E[(X-\mu)^i]$, once again thanks to Maple. From this, in turn, one can get the
{\it scaled moments}, take the limit, and verify {\it asymptotic normality} for as many moments as desired.

We have

{\bf Theorem 3}: 

{\bf (i)} The expected evil-location of a base $b$ random $n$-evil real number is

$$
\frac{2}{b -1}\cdot n \, + \, \frac{b-5}{3 (b-1) } \, +  \, O(\alpha_b^n) \quad,
$$
where $0<\alpha_b<1$.

(When $b=10$ and $n=666$ we get $\frac{2}{9} \cdot 666+\frac{10}{27}=148.1851852\dots$. This differs by $8.054743192\times 10^{-62}$ from the exact value!)

{\bf (ii)} The variance of the  evil-location of a base $b$ random $n$-evil real number is
$$
\frac{2 (b+1)}{3\,(b -1)^2} \cdot n \,+\,
\frac{2 b^2- 14 b + 2 }{9\,(b -1)^2}
+  \, O(\alpha_b^n) \quad .
$$

{\bf (iii)} The third central moment of the  evil-location of a base $b$ random $n$-evil real number is
$$
\frac{2 \left(b +1 \right)^{2}}{3 \left(b -1\right)^{3}} \cdot n + \frac{2 \left(b +1 \right) \left(11 b^{2}-95 b +11\right)}{135 \left(b -1\right)^{3}}
+ O(\alpha_b^n) \, .
$$

It follows that the {\it skewness} (scaled central third moment) tends to $0$ as $n$ goes to infinity.

{\bf (iv)} The fourth central moment of the  evil-location of a base $b$ random $n$-evil real number is
$$
\frac{4 \left(b+1 \right)^{2}}{3 \left(b -1\right)^{4}} \cdot n^2
+
\frac{2 \left(b +1\right) \left(13 b^{2}-30 b +13\right)}{15 \left(b -1\right)^{4}} \cdot n
+
\frac{2 \left(17 b^{2}-11 b +17\right) \left(b^{2}-10 b +1\right)}{135 \left(b -1\right)^{4}}
+ O(\alpha_b^n) \quad .
$$

It follows that the {\it kurtosis} (scaled central fourth moment) tends to $3$ as $n$ goes to infinity.

{\bf (v)} The fifth central moment of  the  evil-location of a base $b$ random $n$-evil real number is
$$
\frac{40 \left(b +1\right)^{3} }{9 \left(b -1\right)^{5}} \cdot n^2
+
\frac{10 \left(31 b^{2}-118 b +31\right) \left(b +1\right)^{2} }{81 \left(b -1\right)^{5}} \cdot n
+
\frac{2 \left(b+1 \right) \left(293 b^{4}-4291 b^{3}+7356 b^{2}-4291 b +293\right)}{1701 \left(b -1\right)^{5}}
+
O(\alpha_b^n) \quad .
$$

It follows that the {\it scaled central fifth moment} tends to $0$ as $n$ goes to infinity.

For the remaining central moments up to the $16$-th see the output file:

{\tt http://sites.math.rutgers.edu/\~{}zeilberg/tokhniot/oHodes6.txt} \quad .

We are happy to report that the limits of the scaled moments approach those of the Gaussian distribution, hence `location of evil-location' over random evil numbers (in any base!) is {\it asymptotically normal}
(at least up to the $16$-th moment).

{\bf Lots and  Lots of Evil Numbers}

It is universally believed, but no one has any clue how to prove it, that every {\it naturally} occurring irrational real number, like $\sqrt{2}$, $\pi$, $e$, and of course, the golden ratio, are {\it normal}.
This means that they behave, statistically, like random real numbers.
The only currently {\it proved} normal real numbers are the (base-dependent) Champernowne's constant, and the non-computable Chaitin constants (see [W]).

Assuming the normality hypothesis, we should expect that the percentage of evil numbers among the first $100000$ prime multiples of $\pi$ should be close to $\%20$, and the
expected evil-location close to $148.185$. This turned out to be true. The output file \hfill\break
{\tt https://sites.math.rutgers.edu/\~{}zeilberg/tokhniot/oHodes7b.txt},
lists which primes $p$, amongst the first $100000$, (from $2$ to $1299709$), are such that $p \cdot \pi$ is evil, along with their respective evil locations.
We found out that $.2010500000$ of them are evil, and the average of their evil-locations is $148.6589406$.

For numbers of the form $\pi \sqrt{x}$, for $ 1 \leq x \leq 10000$ the
ratio is $.2045000000$, see the output file:

{\tt https://sites.math.rutgers.edu/\~{}zeilberg/tokhniot/oHodes2.txt} \quad.

We also explored other bases, see the front of this article:

{\tt https://sites.math.rutgers.edu/\~{}zeilberg/mamarim/mamarimhtml/hodes.html} \quad,

for a few more output files, but you are welcome to explore  your own favorite numbers, once you downloaded {\tt Hodes.txt} \quad .

{\bf References}

[H] Isaac Hodes, {\it The golden ratio is evil}, posted Aug. 24, 2010 \hfill \break
{\tt https://isaachodes.io/p/the-golden-ratio-is-evil/} \quad.

[W] Wikipedia, {\it Normal Number}, {\tt https://en.wikipedia.org/wiki/Normal$\_$number} \quad.

\bigskip
\hrule
\bigskip
Doron Zeilberger, Department of Mathematics, Rutgers University (New Brunswick), Hill Center-Busch Campus, 110 Frelinghuysen
Rd., Piscataway, NJ 08854-8019, USA. {\tt DoronZeil@gmail.com} \hfill\break
{\bf Oct. 4, 2025} 

\end